\documentclass[twoside]{article}%
\usepackage{amssymb}
\usepackage{amsfonts}
\usepackage{amsmath}
\usepackage{graphicx}%
\providecommand{\U}[1]{\protect\rule{.1in}{.1in}}
\newtheorem{theorem}{Theorem}

\newtheorem{corollary}[theorem]{Corollary}

\newtheorem{definition}[theorem]{Definition}

\newtheorem{lemma}[theorem]{Lemma}

\newtheorem{remark}[theorem]{Remark}

\newenvironment{proof}[1][Proof]{\noindent\textbf{#1.} }{\ \rule{0.5em}{0.5em}}
\begin{document}

\title{\textbf{Analytical} \textbf{Blowup Solutions to the }$4$\textbf{-dimensional
Pressureless Navier-Stokes-Poisson Equations with Density-dependent Viscosity}}
\author{Y\textsc{uen} M\textsc{anwai\thanks{E-mail address: nevetsyuen@hotmail.com }}\\\textit{Department of Applied Mathematics, }\\\textit{The Hong Kong Polytechnic University,}\\\textit{Hung Hom, Kowloon, Hong Kong}}
\date{Revised 09-Nov-2008}
\maketitle

\begin{abstract}
We study the 4-dimensional pressureless Navier--Stokes-Poisson equations with
density-dependent viscosity. The analytical solutions with arbitrary time
blowup, in radial symmetry, are constructed in this paper.

\end{abstract}

\section{Introduction}

The evolution of a self-gravitating fluid can be formulated by the
Navier-Stokes-Poisson equations of the following form:
\begin{equation}
\left\{
\begin{array}
[c]{rl}%
{\normalsize \rho}_{t}{\normalsize +\nabla\bullet(\rho u)} & {\normalsize =}%
{\normalsize 0,}\\
{\normalsize (\rho u)}_{t}{\normalsize +\nabla\bullet(\rho u\otimes u)+\nabla
P} & {\normalsize =}{\normalsize -\rho\nabla\Phi+vis(\rho,u),}\\
{\normalsize \Delta\Phi(t,x)} & {\normalsize =\alpha(N)}{\normalsize \rho,}%
\end{array}
\right.  \label{Euler-Poisson}%
\end{equation}
where $\alpha(N)$ is a constant related to the unit ball in $R^{N}$:
$\alpha(1)=2$; $\alpha(2)=2\pi$ and For $N\geq3,$%
\begin{equation}
\alpha(N)=N(N-2)V(N)=N(N-2)\frac{\pi^{N/2}}{\Gamma(N/2+1)},
\end{equation}
where $V(N)$ is the volume of the unit ball in $R^{N}$ and $\Gamma$ is a Gamma
function. And as usual, $\rho=\rho(t,x)$ and $u=u(t,x)\in\mathbf{R}^{N}$ are
the density, the velocity respectively. $P=P(\rho)$\ is the pressure.

In the above system, the self-gravitational potential field $\Phi=\Phi
(t,x)$\ is determined by the density $\rho$ through the Poisson equation.

And $vis(\rho,u)$ is the viscosity function:%
\begin{equation}
vis(\rho,u)=\bigtriangledown(\mu(\rho)\bigtriangledown\bullet u).
\end{equation}
Here we under a common assumption for:
\begin{equation}
\mu(\rho)\doteq\kappa\rho^{\theta}%
\end{equation}
and $\kappa$ and $\theta\geq0$ are the constants. In particular, when
$\theta=0$, it returns the expression for the $u$ dependent only viscosity
function:%
\begin{equation}
vis(\rho,u)=\kappa\Delta u.
\end{equation}
And the vector Laplacian in $u(t,r)$ can be expressed:%
\begin{equation}
\Delta u=u_{rr}+\frac{N-1}{r}u_{r}-\frac{N-1}{r^{2}}u.
\end{equation}
The equations (\ref{Euler-Poisson})$_{1}$ and (\ref{Euler-Poisson})$_{2}$
$(vis(\rho,u)\neq0)$ are the compressible Navier-Stokes equations with forcing
term. The equation (\ref{Euler-Poisson})$_{3}$ is the Poisson equation through
which the gravitational potential is determined by the density distribution of
the density itself. Thus, we call the system (\ref{Euler-Poisson}) the
Navier--Stokes-Poisson equations.

Here, if the $vis(\rho,u)=0$, the system is called the Euler-Poisson
equations.\ In this case, the equations can be viewed as a prefect gas model.
For $N=3$, (\ref{Euler-Poisson}) is a classical (nonrelativistic) description
of a galaxy, in astrophysics. See \cite{C}, \cite{M1} for a detail about the system.

$P=P(\rho)$\ is the pressure. The $\gamma$-law can be applied on the pressure
$P(\rho)$, i.e.%
\begin{equation}
{\normalsize P}\left(  \rho\right)  {\normalsize =K\rho}^{\gamma}\doteq
\frac{{\normalsize \rho}^{\gamma}}{\gamma}, \label{gamma}%
\end{equation}
which is a commonly the hypothesis. The constant $\gamma=c_{P}/c_{v}\geq1$,
where $c_{P}$, $c_{v}$\ are the specific heats per unit mass under constant
pressure and constant volume respectively, is the ratio of the specific heats,
that is, the adiabatic exponent in (\ref{gamma}). In particular, the fluid is
called isothermal if $\gamma=1$. With $K=0$, we call the system is pressureless.

For the $3$-dimensional case, we are interested in the hydrostatic equilibrium
specified by $u=0$. According to \cite{C}, the ratio between the core density
$\rho(0)$ and the mean density $\overset{\_}{\rho}$ for $6/5<\gamma<2$\ is
given by%
\begin{equation}
\frac{\overset{\_}{\rho}}{\rho(0)}=\left(  \frac{-3}{z}\dot{y}\left(
z\right)  \right)  _{z=z_{0}}%
\end{equation}
where $y$\ is the solution of the Lane-Emden equation with $n=1/(\gamma-1)$,%
\begin{equation}
\ddot{y}(z)+\dfrac{2}{z}\dot{y}(z)+y(z)^{n}=0,\text{ }y(0)=\alpha>0,\text{
}\dot{y}(0)=0,\text{ }n=\frac{1}{\gamma-1},
\end{equation}
and $z_{0}$\ is the first zero of $y(z_{0})=0$. We can solve the Lane-Emden
equation analytically for%
\begin{equation}
y_{anal}(z)\doteq\left\{
\begin{array}
[c]{ll}%
1-\frac{1}{6}z^{2}, & n=0;\\
\dfrac{\sin z}{z}, & n=1;\\
\dfrac{1}{\sqrt{1+z^{2}/3}}, & n=5,
\end{array}
\right.
\end{equation}
and for the other values, only numerical values can be obtained. It can be
shown that for $n<5$, the radius of polytropic models is finite; for $n\geq5$,
the radius is infinite.

Gambin \cite{G} and Bezard \cite{B} obtained the existence results about the
explicitly stationary solution $\left(  u=0\right)  $ for $\gamma=6/5$ in
Euler-Poisson equations$:$%
\begin{equation}
\rho=\left(  \frac{3KA^{2}}{2\pi}\right)  ^{5/4}\left(  1+A^{2}r^{2}\right)
^{-5/2}, \label{stationsoluionr=6/5}%
\end{equation}
where $A$ is constant.\newline The Poisson equation (\ref{Euler-Poisson}%
)$_{3}$ can be solved as%
\begin{equation}
{\normalsize \Phi(t,x)=}\int_{R^{N}}G(x-y)\rho(t,y){\normalsize dy,}%
\end{equation}
where $G$ is the Green's function for the Poisson equation in the
$N$-dimensional spaces defined by
\begin{equation}
G(x)\doteq\left\{
\begin{array}
[c]{ll}%
|x|, & N=1;\\
\log|x|, & N=2;\\
\dfrac{-1}{|x|^{N-2}}, & N\geq3.
\end{array}
\right.
\end{equation}
In the following, we always seek solutions in radial symmetry. Thus, the
Poisson equation (\ref{Euler-Poisson})$_{3}$ is transformed to%
\begin{equation}
{\normalsize r^{N-1}\Phi}_{rr}\left(  {\normalsize t,x}\right)  +\left(
N-1\right)  r^{N-2}\Phi_{r}{\normalsize =}\alpha\left(  N\right)
{\normalsize \rho r^{N-1},}%
\end{equation}%
\begin{equation}
\Phi_{r}=\frac{\alpha\left(  N\right)  }{r^{N-1}}\int_{0}^{r}\rho
(t,s)s^{N-1}ds.
\end{equation}

\begin{definition}
[Blowup]We say a solution blows up if one of the following conditions is
satisfied:\newline(1)The solution becomes infinitely large at some point $x$
and some finite time $T_{0}$;\newline(2)The derivative of the solution becomes
infinitely large at some point $x$ and some finite time $T_{0}$.
\end{definition}

In this paper, we concern about blowup solutions for the $4$-dimensional
pressureless Navier--Stokes-Poisson equations with the density-dependent
viscosity. And our aim is to construct a family of such blowup solutions.

Historically in astrophysics, Goldreich and Weber \cite{GW} constructed the
analytical blowup solution (collapsing) of the $3$-dimensional Euler-Poisson
equations for $\gamma=4/3$ for the non-rotating gas spheres. After that,
Makino \cite{M1} obtained the rigorously mathematical proof of the existence
of such kind of blowup solutions. And in \cite{DXY}, we find the extension of
the above blowup solutions to the case . In \cite{Y}, the solutions with a
from is rewritten as

For $N\geq3$ and $\gamma=(2N-2)/N$,
\begin{equation}
\left\{
\begin{array}
[c]{c}%
\rho(t,r)=\left\{
\begin{array}
[c]{c}%
\dfrac{1}{a(t)^{N}}y(\frac{r}{a(t)})^{N/(N-2)},\text{ for }r<a(t)Z_{\mu};\\
0,\text{ for }a(t)Z_{\mu}\leq r.
\end{array}
\right.  \text{, }{\normalsize u(t,r)=}\dfrac{\dot{a}(t)}{a(t)}%
{\normalsize r,}\\
\ddot{a}(t){\normalsize =}-\dfrac{\lambda}{a(t)^{N-1}},\text{ }%
{\normalsize a(0)=a}_{0}>0{\normalsize ,}\text{ }\dot{a}(0){\normalsize =a}%
_{1},\\
\ddot{y}(z){\normalsize +}\dfrac{N-1}{z}\dot{y}(z){\normalsize +}\dfrac
{\alpha(N)}{(2N-2)K}{\normalsize y(z)}^{N/(N-2)}{\normalsize =\mu,}\text{
}y(0)=\alpha>0,\text{ }\dot{y}(0)=0,
\end{array}
\right.  \label{solution2}%
\end{equation}
where $\mu=[N(N-2)\lambda]/(2N-2)K$ and the finite $Z_{\mu}$ is the first zero
of $y(z)$;

For $N=2$ and $\gamma=1,$%
\begin{equation}
\left\{
\begin{array}
[c]{c}%
\rho(t,r)=\dfrac{1}{a(t)^{2}}e^{y(r/a(t))}\text{, }{\normalsize u(t,r)=}%
\dfrac{\dot{a}(t)}{a(t)}{\normalsize r;}\\
\ddot{a}(t){\normalsize =}-\dfrac{\lambda}{a(t)},\text{ }{\normalsize a(0)=a}%
_{0}>0{\normalsize ,}\text{ }\dot{a}(0){\normalsize =a}_{1};\\
\ddot{y}(z){\normalsize +}\dfrac{1}{z}\dot{y}(z){\normalsize +\dfrac{2\pi}%
{K}e}^{y(z)}{\normalsize =\mu,}\text{ }y(0)=\alpha,\text{ }\dot{y}(0)=0,
\end{array}
\right.  \label{solution 3}%
\end{equation}
where $K>0$, $\mu=2\lambda/K$ with a sufficiently small $\lambda$ and $\alpha$
are constants.

For the construction of special analytical solutions to the Navier-Stokes
equations in $R^{N}$ or the Navier-Stokes-Poisson equations in $R^{3}$ without
pressure with $\theta=1$, readers may refer Yuen's recent results in
\cite{Y2}, \cite{Y3} respectively.

In this article, the analytical blowup solutions are constructed in the
pressureless Euler-Poisson equations with density-dependent viscosity in
$R^{4}$ with $\theta=5/4$ in radial symmetry:%
\begin{equation}
\left\{
\begin{array}
[c]{rl}%
\rho_{t}+u\rho_{r}+\rho u_{r}+{\normalsize \dfrac{3}{r}\rho u} &
{\normalsize =0,}\\
\rho\left(  u_{t}+uu_{r}\right)  & {\normalsize =-}\dfrac{\alpha(4)\rho}%
{r^{3}}%
%TCIMACRO{\dint _{0}^{r}}%
%BeginExpansion
{\displaystyle\int_{0}^{r}}
%EndExpansion
\rho(t,s)s^{3}ds+[\kappa\rho^{5/4}]_{r}u_{r}+(\kappa\rho^{5/4})(u_{rr}%
+\dfrac{3}{r}u_{r}-\dfrac{3}{r^{2}}u),
\end{array}
\right.  \label{gamma=1}%
\end{equation}
in the form of the following theorem.

\begin{theorem}
\label{thm:1}For the $4$-dimensional pressureless Navier--Stokes-Poisson
equations with $\theta=5/4$, in radial symmetry, (\ref{gamma=1}), there exists
a family of blowup solutions,%
\begin{equation}
\left\{
\begin{array}
[c]{c}%
\rho(t,r)=\dfrac{1}{(T-Ct)^{4}}y(\frac{r}{T-Ct})^{4},\text{ }%
{\normalsize u(t,r)=}\dfrac{-C}{T-Ct}{\normalsize r;}\\
\ddot{y}(z){\normalsize +}\dfrac{3}{z}\dot{y}(z){\normalsize +\dfrac
{\alpha(4)}{\kappa C}y(z)}^{4}{\normalsize =0,}\text{ }y(0)=\alpha,\text{
}\dot{y}(0)=0,
\end{array}
\right.  \label{solution1}%
\end{equation}
where $T>0$, $\kappa>0$, $C>0$ and $\alpha$ are constants.\newline And the
solutions blow up in the finite time $T/C$.
\end{theorem}

\section{Separable Blowup Solutions}

Before presenting the proof of Theorem \ref{thm:1}, we prepare some lemmas.
First, we obtain the solutions for the continuity equation of mass in radial
symmetry (\ref{gamma=1})$_{1}$.

\begin{lemma}
\label{lem:generalsolutionformasseq}For the 4-dimensional conservation of mass
in radial symmetry
\begin{equation}
\rho_{t}+u\rho_{r}+\rho u_{r}+\dfrac{3}{r}\rho u=0,
\label{massequationspherical}%
\end{equation}
there exist solutions,%
\begin{equation}
\rho(t,r)=\frac{1}{(T-Ct)^{4}}y(\frac{r}{T-Ct})^{4},\text{ }%
{\normalsize u(t,r)=}\frac{-C}{T-Ct}{\normalsize r,}
\label{generalsolutionformassequation}%
\end{equation}
where $T$ and $C$ are positive constants.
\end{lemma}

\begin{proof}
We just plug (\ref{generalsolutionformassequation}) into
(\ref{massequationspherical}). Then
\begin{align*}
&  \rho_{t}+\rho_{r}u+\rho u_{r}+\dfrac{3}{r}\rho u\\
&  =\frac{(-4)(-C)y(\frac{r}{T-Ct})^{4}}{(T-Ct)^{4}}+\frac{4y(\frac{r}%
{T-Ct})^{3}\dot{y}(\frac{r}{T-Ct})}{(T-Ct)^{4}}\frac{r(-1)(-C)}{(T-Ct)^{2}}\\
&  +\frac{4y(\frac{r}{T-Ct})^{3}\dot{y}(\frac{r}{T-Ct})}{(T-Ct)^{4}}\frac
{1}{T-Ct}\frac{(-C)}{T-Ct}{\normalsize r}+\frac{y(\frac{r}{T-Ct})^{4}%
}{(T-Ct)^{4}}\frac{(-C)}{T-Ct}{\normalsize +}\dfrac{3}{r}\frac{y(\frac
{r}{T-Ct})^{4}}{(T-Ct)^{3}}\frac{(-C)}{T-Ct}{\normalsize r}\\
&  =\frac{4Cy(\frac{r}{T-Ct})^{4}}{(T-Ct)^{4}}+\frac{4Cy(\frac{r}{T-Ct}%
)^{3}\dot{y}(\frac{r}{T-Ct})}{(T-Ct)^{6}}r\\
&  -\frac{4Cy(\frac{r}{T-Ct})^{3}\dot{y}(r/(T-Ct))}{(T-Ct)^{6}}r-\frac
{Cy(\frac{r}{T-Ct})^{4}}{(T-Ct)^{4}}-\frac{3Cy(\frac{r}{T-Ct})^{4}}%
{(T-Ct)^{4}}\\
&  =0.
\end{align*}
The proof is completed.
\end{proof}

Besides, we need the lemma for stating the property of the function $y(z)$.
The similar lemma was already given in Lemmas 2.1, \cite{DXY}, by the fixed
point theorem. The proof is similar, the proof may be skipped here.

\begin{lemma}
\label{lemma2}For the ordinary differential equation,%
\begin{equation}
\left\{
\begin{array}
[c]{c}%
\ddot{y}(z)+\dfrac{3}{z}\overset{\cdot}{y}(z)+\frac{\alpha(4)}{5C\kappa
}y(z)^{4}{\normalsize =0,}\\
y(0)=\alpha>0,\text{ }\dot{y}(0)=0,
\end{array}
\right.  \label{SecondorderElliptic}%
\end{equation}
where $\alpha(4)$, $C$ and $\kappa$ are positive constants,\newline has a
solution $y(z)\in C^{2}$ provided that $y(z)\subset\lbrack\alpha,0]$.
\end{lemma}

Here we are already to give the proof of Theorem \ref{thm:1}.

\begin{proof}
[Proof of Theorem 2]From Lemma \ref{lem:generalsolutionformasseq}, it is clear
for that (\ref{solution1}) satisfy (\ref{gamma=1})$_{1}$. For the momentum
equation (\ref{gamma=1})$_{2}$, we get,%
\begin{align}
&  \rho(u_{t}+uu_{r})+\frac{\alpha(4)\rho}{r^{3}}%
%TCIMACRO{\dint \limits_{0}^{r}}%
%BeginExpansion
{\displaystyle\int\limits_{0}^{r}}
%EndExpansion
\rho(t,s)s^{3}ds-[\kappa\rho^{5/4}]_{r}u_{r}-\kappa\rho^{5/4}(u_{rr}+\dfrac
{3}{r}u_{r}-\dfrac{3}{r^{2}}u)\\
&  =\rho\left[  \frac{(-C)(-1)(-C)}{(T-Ct)^{2}}r+\frac{(-C)}{T-Ct}r\cdot
\frac{(-C)}{T-Ct}\right]  +\frac{\alpha(4)\rho}{r^{3}}%
%TCIMACRO{\dint \limits_{0}^{r}}%
%BeginExpansion
{\displaystyle\int\limits_{0}^{r}}
%EndExpansion
\frac{y(\frac{s}{T-Ct})^{4}}{(T-Ct)^{4}}s^{3}ds\\
&  -\left[  \kappa\frac{1}{(T-Ct)^{4}}y\left(  \frac{r}{T-Ct}\right)
^{4}\right]  _{r}^{5/4}\frac{(-C)}{T-Ct}-0\nonumber\\
&  =\frac{\alpha(4)\rho}{r^{3}}%
%TCIMACRO{\dint \limits_{0}^{r}}%
%BeginExpansion
{\displaystyle\int\limits_{0}^{r}}
%EndExpansion
\frac{y(\frac{s}{T-ct})^{4}}{(T-Ct)^{4}}s^{3}ds-\frac{5}{4}\kappa\left[
\frac{1}{(T-Ct)^{4}}y\left(  \frac{r}{T-Ct}\right)  ^{4}\right]  ^{1/4}%
\frac{4y\left(  \frac{r}{T-Ct}\right)  ^{3}\overset{\cdot}{y}\left(  \frac
{r}{T-Ct}\right)  }{(T-Ct)^{4}}\frac{1}{T-Ct}\frac{(-C)}{T-Ct}\\
&  =\frac{\alpha(4)\rho}{r^{3}}%
%TCIMACRO{\dint \limits_{0}^{r}}%
%BeginExpansion
{\displaystyle\int\limits_{0}^{r}}
%EndExpansion
\frac{y(\frac{s}{T-Ct})^{4}}{(T-Ct)^{4}}s^{3}ds+5C\kappa\frac{y\left(
\frac{r}{T-Ct}\right)  }{T-Ct}\frac{y\left(  \frac{r}{T-Ct}\right)
^{3}\overset{\cdot}{y}\left(  \frac{r}{T-Ct}\right)  }{(T-Ct)^{4}}\frac
{1}{(T-Ct)^{2}}\\
&  =\frac{\alpha(4)\rho}{r^{3}}%
%TCIMACRO{\dint \limits_{0}^{r}}%
%BeginExpansion
{\displaystyle\int\limits_{0}^{r}}
%EndExpansion
\frac{y(\frac{s}{T-Ct})^{4}}{(T-Ct)^{4}}s^{3}ds+\frac{5C\kappa\rho}%
{(T-Ct)^{3}}\overset{\cdot}{y}\left(  \frac{r}{T-Ct}\right) \\
&  =\frac{\rho}{(T-Ct)^{3}}\left[  5C\kappa\overset{\cdot}{y}(\frac{r}%
{T-Ct})+\frac{\alpha(4)}{r^{3}(T-Ct)}%
%TCIMACRO{\dint \limits_{0}^{r}}%
%BeginExpansion
{\displaystyle\int\limits_{0}^{r}}
%EndExpansion
y(\frac{s}{T-Ct})^{4}s^{3}ds\right] \\
&  =\frac{\rho}{(T-Ct)^{3}}\left[  5C\kappa\overset{\cdot}{y}(\frac{r}%
{T-Ct})+\frac{\alpha(4)}{(\frac{r}{T-Ct})^{3}}%
%TCIMACRO{\dint \limits_{0}^{r/(T-Ct)}}%
%BeginExpansion
{\displaystyle\int\limits_{0}^{r/(T-Ct)}}
%EndExpansion
y(s)^{4}s^{3}ds\right] \\
&  =\frac{\rho}{(T-Ct)^{3}}Q\left(  \frac{r}{T-Ct}\right)  .
\end{align}
And denote%
\begin{equation}
Q(\frac{r}{T-Ct})\doteq{\normalsize Q(z)=}5C\kappa\overset{\cdot}%
{y}(z){\normalsize +}\frac{\alpha(4)}{z^{3}}%
%TCIMACRO{\dint \limits_{0}^{z}}%
%BeginExpansion
{\displaystyle\int\limits_{0}^{z}}
%EndExpansion
y{\normalsize (s)}^{4}{\normalsize s}^{3}{\normalsize ds.}%
\end{equation}
Differentiate $Q(z)$\ with respect to $z$,%
\begin{align}
&  \overset{\cdot}{Q}(z)\\
&  =5C\kappa\overset{\cdot\cdot}{y}(z){\normalsize +}\alpha
(4)y{\normalsize (z)}^{4}+\frac{(-3)\alpha(4)}{z^{4}}%
%TCIMACRO{\dint \limits_{0}^{z}}%
%BeginExpansion
{\displaystyle\int\limits_{0}^{z}}
%EndExpansion
y{\normalsize (s)}^{4}{\normalsize s}^{3}{\normalsize ds}\\
&  =-\frac{3}{z}\cdot5C\kappa{\normalsize \overset{\cdot}{y}(z)+}\frac
{(-3)}{z}\cdot\frac{\alpha(4)}{z^{3}}%
%TCIMACRO{\dint \limits_{0}^{z}}%
%BeginExpansion
{\displaystyle\int\limits_{0}^{z}}
%EndExpansion
y{\normalsize (s)}^{4}{\normalsize s}^{3}{\normalsize ds}\\
&  =\frac{-3}{z}Q(z),
\end{align}
where the above result is due to the fact that we choose the following
ordinary differential equation,%
\begin{equation}
\left\{
\begin{array}
[c]{c}%
\ddot{y}(z)+\dfrac{3}{z}\overset{\cdot}{y}(z)+\dfrac{\alpha(4)}{5C\kappa
}y(z)^{4}=0.\\
{\normalsize y(0)=\alpha>0,}\text{ }\dot{y}(0){\normalsize =0.}%
\end{array}
\right.
\end{equation}
With $Q(0)=0$, this implies that $Q(z)=0$. Thus, the momentum equation
(\ref{gamma=1})$_{2}$ is satisfied.\newline Now we are able to show that the
family of the solutions blow up, in the finite time $T/C$. This completes the proof.
\end{proof}

The statement about the blowup rate will be immediately followed:

\begin{corollary}
The blowup rate of the solution (\ref{solution1}) is
\begin{equation}
\underset{t\rightarrow T/C^{-}}{\lim}\rho(t,0)(T-Ct)^{4}\geq O(1).
\end{equation}

\end{corollary}

Given that the sign of the constant $C$ in (\ref{solution 3}) is changed to be
negative, the below corollary is clearly shown.

\begin{corollary}
For the $4$-dimensional pressureless Navier--Stokes-Poisson equations, with
$\theta=5/4$, in radial symmetry, (\ref{gamma=1}), there exists a family of
solutions,%
\begin{equation}
\left\{
\begin{array}
[c]{c}%
\rho(t,r)=\dfrac{1}{(T-Ct)^{4}}y\left(  \frac{r}{T-Ct}\right)  ^{4}\text{,
}{\normalsize u(t,r)=}\dfrac{-C}{T-Ct}{\normalsize r;}\\
\ddot{y}(z){\normalsize +}\dfrac{3}{z}\dot{y}(z){\normalsize +\dfrac
{\alpha(4)}{5C\kappa}y(z)^{4}=0,}\text{ }y(0)=\alpha>0,\text{ }\dot{y}(0)=0,
\end{array}
\right.
\end{equation}
where $T>0$, $\kappa>0$, $C<0$ and $\alpha$ are constants.
\end{corollary}

\begin{remark}
Besides, if we consider the $4$-dimensional Navier-Stokes equations with the
repulsive force in radial symmetry,%
\begin{equation}
\left\{
\begin{array}
[c]{rl}%
\rho_{t}+u\rho_{r}+\rho u_{r}+{\normalsize \dfrac{3}{r}\rho u} &
{\normalsize =0,}\\
\rho\left(  u_{t}+uu_{r}\right)  & {\normalsize =+}\dfrac{\alpha(4)\rho}%
{r^{3}}\int_{0}^{r}\rho(t,s)s^{3}ds+[\kappa\rho^{5/4}]u_{r}+(\kappa\rho
^{5/4})(u_{rr}+\dfrac{3}{r}u_{r}-\dfrac{3}{r^{2}}u),
\end{array}
\right.
\end{equation}
the special solutions are:%
\begin{equation}
\left\{
\begin{array}
[c]{c}%
\rho(t,r)=\dfrac{1}{(T-Ct)^{4}}y\left(  \frac{r}{T-Ct}\right)  ^{4}\text{,
}{\normalsize u(t,r)=}\dfrac{-C}{T-Ct}{\normalsize r}\\
\ddot{y}(z){\normalsize +}\dfrac{3}{z}\dot{y}(z)-{\normalsize \frac{\alpha
(4)}{5C\kappa}}y\left(  z\right)  ^{4}{\normalsize =0,}\text{ }y(0)=\alpha
,\text{ }\dot{y}(0)=0.
\end{array}
\right.
\end{equation}

\end{remark}

\end{document}